\documentclass[graybox]{svmult}
\pdfoutput=1

\usepackage{type1cm}        
\usepackage{makeidx}         
\usepackage{graphicx}        
\usepackage{multicol}        
\usepackage[bottom]{footmisc}

\usepackage{newtxtext}       %
\usepackage{newtxmath}       

\usepackage{mathrsfs}
\usepackage[capitalise]{cleveref}
\usepackage[labelfont=bf]{caption}
\usepackage{subfig}

\usepackage{tikz}
\usepackage{pgfplots}
\usepackage{xcolor}

\usepackage{enumerate}
\usepackage{algorithm, algorithmic}
	
\definecolor{MPIblue}{RGB}{51,165,195}
\colorlet{mpiblue}{MPIblue}
\definecolor{MPIgrey}{RGB}{135,135,141}
\colorlet{mpigrey}{MPIgrey}
\definecolor{MPIgreen}{RGB}{0,118,117}
\colorlet{mpigreen}{MPIgreen}
\definecolor{MPIred}{RGB}{120,0,75}
\colorlet{mpired}{MPIred}
\definecolor{MPIsand}{RGB}{236,233,212}
\colorlet{mpisand}{MPIsand}
\definecolor{MPItext}{RGB}{56,60,60}
\colorlet{mpitext}{MPItext}
\definecolor{MPIbluetext}{RGB}{23,161,193}
\colorlet{mpibluetext}{MPIbluetext}
\definecolor{MPIbackground}{RGB}{239,238,227}
\colorlet{alg1col1}{MPIblue!30}
\colorlet{alg1col2}{MPIsand}

\newcommand{\pmat}[1]{\ensuremath{\textbf{{#1}}(\mu)}}
\newcommand{\pmatr}[1]{\ensuremath{\textbf{{#1}}_{r}(\mu)}}
\newcommand\Tstrut{\rule{0pt}{2.6ex}}         

\makeindex             

\begin{document}

\title*{An Adaptive Sampling Approach for the Reduced Basis Method}
\author{Sridhar Chellappa, Lihong Feng and Peter Benner}
\institute{Sridhar Chellappa \at Max Planck Institute for Dynamics of Complex Technical Systems, Magdeburg, Germany. \email{chellappa@mpi-magdeburg.mpg.de}.\\ Supported by the International Max Planck Research School for Advanced Methods in Process and Systems Engineering (IMPRS-ProEng).
\and Lihong Feng \at Max Planck Institute for Dynamics of Complex Technical Systems, Magdeburg,  Germany. \email{feng@mpi-magdeburg.mpg.de}.
\and Peter Benner \at Max Planck Institute for Dynamics of Complex Technical Systems, Magdeburg,  Germany and
Fakult{\"a}t f{\"u}r Mathematik, Otto-von-Guericke Universit{\"a}t Magdeburg, Germany. \email{benner@mpi-magdeburg.mpg.de}.}
%
%
\maketitle

\abstract*{The offline time of the reduced basis method can be very long given a large training set of parameter samples. This usually happens when the system has more than two independent parameters. On the other hand, if the training set includes fewer parameter samples, the greedy algorithm might produce a reduced-order model with large errors at the samples outside of the training set. We introduce a method based on a surrogate error model to efficiently sample the parameter domain such that the training set is adaptively updated starting from a coarse set with a small number of parameter samples. A sharp a posteriori error estimator is evaluated on a coarse training set. Radial Basis Functions are used to interpolate the error estimator over a separate fine training set. Points from the fine training set are added into the coarse training set at every iteration based on a user defined criterion. In parallel, parameter samples satisfying a defined tolerance are adaptively removed from the coarse training set. The approach is shown to avoid high computational costs by using a small training set and to provide a reduced-order model with guaranteed accuracy over a fine training set. Further, we show numerical evidence that the reduced-order model meets the defined tolerance over an independently sampled test set from the parameter domain.}

\abstract{The offline time of the reduced basis method can be very long given a large training set of parameter samples. This usually happens when the system has more than two independent parameters. On the other hand, if the training set includes fewer parameter samples, the greedy algorithm might produce a reduced-order model with large errors at the samples outside of the training set. We introduce a method based on a surrogate error model to efficiently sample the parameter domain such that the training set is adaptively updated starting from a coarse set with a small number of parameter samples. A sharp a posteriori error estimator is evaluated on a coarse training set. Radial Basis Functions are used to interpolate the error estimator over a separate fine training set. Points from the fine training set are added into the coarse training set at every iteration based on a user defined criterion. In parallel, parameter samples satisfying a defined tolerance are adaptively removed from the coarse training set. The approach is shown to avoid high computational costs by using a small training set and to provide a reduced-order model with guaranteed accuracy over a fine training set. Further, we show numerical evidence that the reduced-order model meets the defined tolerance over an independently sampled test set from the parameter domain.}
\section{Introduction}
\label{sec:1}
The need of modeling physical processes accurately very often leads to mathematical models of large-scale. 
Simulating such large-scale systems poses challenges on computer memory and computational power. Model Order Reduction (MOR) aims to speedup simulation of the large-scale full order model (FOM) by constructing a small-scale system, called the reduced-order model (ROM), which serves as a good approximation of the FOM. 

For FOMs arising from discretization of parametric partial differential equations (PDEs), the Reduced Basis (RB) method is  a popular choice \cite{morQuaMN16}. It is suitable especially when the FOM needs to be solved repeatedly for a range of parameters. The RB method involves two stages. The first is an expensive \emph{offline stage} in which a subspace of the solution manifold is constructed by simulations of the FOM at certain chosen parameters. Then, using Galerkin projection, the solution to the PDE is projected onto this subspace to generate the ROM. In the second \emph{online stage}, instead of solving the FOM, one solves the ROM for any parameter of interest. The RB method works well when the solution manifold is supposed to have a small Kolmogorov $n$-width. In this work, we focus on the offline stage of the RB method.

The offline stage in the RB method requires the construction of a representative set of parameters from the domain of interest. This is called the \emph{training set}. The subspace approximating the solution manifold is constructed by solving the FOM at a set of selected points from this training set, picked using a greedy algorithm \cite{morHaaO08}. The choice of a training set is nontrivial. On the one hand, if it includes too few parameters, the original solution manifold may not be adequately represented, leading to a poor ROM with large error. On the other hand, if it is too fine, the offline time can be prohibitively long. When the PDE involves several parameters, properly defining the training set can be a severe computational issue. 

Several authors have attempted to address this issue. The earliest work to consider an adaptive sampling of the training set was \cite{morSen08}, where the author proposes a multi-stage method. The greedy algorithm is run several times over randomly sampled small training sets to generate the RB basis. Then, the ROM is tested over a much larger training set and the greedy algorithm is re-run only on those points failing the tolerance criterion.
The authors of \cite{morHaaDO11} address the issue of large training sets by two approaches. The first one is a procedure to monitor the error over an additional validation parameter set. If a large error is detected, then the training set is further refined, either uniformly or locally. The second approach, similar to the one presented in \cite{morEftKP11}, is based on partitioning the parameter domain adaptively, and generating  a local basis for each partition. Other approaches for adaptive sampling are proposed in \cite{morHesSZ14} and \cite{morMadS13}. 
More recently, an approach based on Kriging interpolation and clustering is proposed in \cite{morTaineA15}, to tackle the problem of high-dimensional parameter spaces. An interpolant of the residual norm is calculated over a fine grid of parameters. Then, k-means clustering is used to identify parameters that have high probability of presenting larger errors.   

The approach we propose here is based on interpolating a sharp \emph{a posteriori} output error estimator originally proposed in \cite{morCheFB19a}. The RB method is initialized with a coarse training set. We update this set progressively by adding or removing points from it. At each iteration of the greedy algorithm, we estimate the error of the ROM at every parameter in the coarse training set. Note that the evaluation of the estimator does not need to evaluate the FOM. In order to further reduce the computational costs, we then interpolate the estimated error over a fine training set and use the interpolant as the error surrogate, which will replace the error estimator for estimating the ROM error over the fine training set. To achieve this, a surrogate model based on Radial Basis Functions (RBFs) is used.

At each iteration of the greedy algorithm, the ROM errors at the parameters in the fine training set are checked by the error surrogate. Those parameters corresponding to large values of the error surrogate are selected and added into the coarse training set. The error surrogate is much cheaper to compute than the error estimator. Therefore, using the former instead of the latter to check the ROM error over the fine training set, will reduce the computational cost. Additionally, if any parameter in the coarse training set achieves the required ROM accuracy, we remove it from the coarse training set. Such an approach is able to construct a small, representative training set by fully exploring the parameter domain with reduced computational cost. Forming the RBF interpolant is a relatively cheap operation with the most expensive contribution coming from the need to solve an $\ell \times \ell$ linear system of equations at each parameter sample, where $\ell$ is the cardinality of the coarse training set, which is expected to be small.  The algorithm proposed in our work is based on the use of the RBF interpolant and offers certain advantages over the most relevant methods in \cite{morHesSZ14} and \cite{morTaineA15}. Our proposed approach differs in the following aspects, when compared to \cite{morHesSZ14}:
	\begin{itemize}
		\item We use a primal-dual a posteriori error estimator on a coarse training set and a cheaply computable error surrogate on the fine training set.  No cheap error estimator is used on the fine training set in~\cite{morHesSZ14}, potentially leading to a larger computational effort.
		\item At each iteration of the greedy algorithm, a saturation assumption is introduced in ~\cite{morHesSZ14} in order to avoid calculation of the ROM at certain parameters. However, this requires estimation of a saturation constant which needs to be defined a priori by the user. In our method, no such constant needs to be estimated.
	\end{itemize}
Furthermore, when compared to \cite{morTaineA15},
	\begin{itemize}
		\item We consider both addition and removal of parameter samples from the coarse training set. The method in \cite{morTaineA15} is restricted to adding new samples only.
		\item The use of Kriging interpolation involves the estimation of several hyper-parameters. This can be a computational bottleneck \cite{Steetal17}. In the case of RBF, there is only one free parameter at most. If using some special kernel functions, e.g. polyharmonic spline kernels, there are no free parameters to tune. 
		\item Finally, a nonlinear model was considered in~\cite{morTaineA15} in order to demonstrate the adaptive sampling approach. An offline hyper-reduction step was used to reduce the complexity of the online nonlinearity evaluations, where a second training set for the nonlinear function is used in the hyper-reduction step, and is likely to be separately given and fixed. In our approach, we employ the Adaptive POD-Greedy-(D)EIM algorithm proposed in \cite{morCheFB19a}. This avoids the need for a potentially separate training set for the hyper-reduction phase. Instead, the (D)EIM basis generation is carried out as a part of the greedy loop. In this way, we propose a fully adaptive approach in~\cref{alg:3}.
	\end{itemize}
The rest of the paper is organized as follows. In~\cref{sec:2}, the basic idea behind the RB method is reviewed. An adaptive parameter sampling approach using the RBF interpolation is proposed and elaborated in~\cref{sec:3}. There, we begin with a brief introduction to the theory of RBFs. Then an adaptive sampling approach for the standard POD-Greedy algorithm is proposed to show a general framework for the proposed technique. In Subsection~\ref{subsec:fulladapt}, a fully adaptive algorithm is proposed, which integrates the adaptive sampling technique with an adaptive basis construction method: the adaptive POD-Greedy-(D)EIM algorithm for general parametric nonlinear/nonaffine systems. A strategy of computing the shape parameter for RBF interpolation is detailed in Subsection~\ref{subsec:RBF_shape}. In \cref{sec:4}, we show numerical results for a series of benchmark examples,  validating the proposed algorithms. We conclude in \cref{sec:5} by summarizing our results and suggesting future research directions. 
\section{Reduced Basis Method}
\label{sec:2}
For the sake of completeness, we briefly highlight the key idea of the RB method. Consider a nonlinear parametric dynamical system arising from space-time discretization of a model of PDEs,
	\begin{equation}
		\begin{aligned}
			\pmat{E} x^{k+1}(\mu) &= \pmat{A} x^{k}(\mu) + f(x^{k}(\mu);\mu) + \pmat{B} u^{k},\\
			y^{k+1}(\mu)  &= \pmat{C} x^{k+1}(\mu)
		\end{aligned}
		\label{eq:fom}
	\end{equation}
where $x^{k}(\mu) \in \mathbb{R}^{n}$ is the state vector at the $k$-th time step $t^{k}$ in the time interval $[0,T]$ divided as $t^{0} < t^{1} < \cdots < t^{k} < \cdots < t^{K} = T$. $\pmat{E},\, \pmat{A} \in \mathbb{R}^{n \times n}$ are the parameter dependent system matrices. $f(x^{k}(\mu);\mu) \in \mathbb{R}^{n}$ depicts the state, parameter dependent system nonlinearity. $\pmat{B} \in \mathbb{R}^{n \times q}, \pmat{C} \in \mathbb{R}^{m \times n}$ are the input and output matrices, respectively, while $u^{k}$ denotes the input at the $k$-th time step
and $y^{k+1}(\mu)$ denotes the quantities-of-interest. The vector of parameters $\mu\in \mathbb{R}^{d}$ belongs to a parameter domain $\mathscr{P} \subset \mathbb{R}^{d}$. The space discretization can be done through any method of choice such as Finite Element Method (FEM), Finite Volume Method (FVM) or the Finite Difference Method (FDM). For the time discretization, we adopt a semi-implicit  scheme for ease of using the a posteriori error estimator in~\cite{morCheFB19a}. The RB method attempts to generate a subspace $\mathscr{X}$ with dimension $r\ll n$, that best approximates the solution manifold. Let $\textbf{V} \subset \mathbb{R}^{n \times r}$ be an orthogonal basis of $\mathscr{X}$. The ROM is then given via Galerkin projection,
	\begin{equation}
		\begin{aligned}
			\pmatr{E} x_{r}^{k+1}(\mu) &= \pmatr{A} x_{r}^{k}(\mu) + f_{r}(x_{r}^{k}(\mu);\mu) + \pmatr{B} u^{k},\\
			y_{r}^{k+1}(\mu)  &= \pmatr{C} x_{r}^{k+1}(\mu),
		\end{aligned}
		\label{eq:rom}
	\end{equation}
where $x_{r}^{k} \in \mathbb{R}^{r}$ is the reduced state vector at the $k$-th time step. $\pmatr{E} \coloneqq \textbf{V}^{T} \pmat{E} \textbf{V} \in \mathbb{R}^{r \times r}$, $\pmatr{A} \coloneqq \textbf{V}^{T} \pmat{A} \textbf{V} \in \mathbb{R}^{r \times r}$ are the reduced system matrices. $\pmatr{B} \coloneqq \textbf{V}^{T} \pmat{B} \in \mathbb{R}^{r \times q}$, $\pmatr{C} \coloneqq \pmat{C} \textbf{V} \in \mathbb{R}^{m \times r}$ are the reduced input and output matrices, respectively. The reduced nonlinear term is given by $f_{r}(\textbf{V} x_{r}^{k}(\mu); \mu ) \coloneqq \textbf{V}^{T} f(\textbf{V} x_{r}^{k}(\mu); \mu ) \in \mathbb{R}^{r}$. An important assumption is that the system matrices and the input, output matrices have affine decomposition, i.e., any of them can be formulated as a sum of products of parameter dependent and parameter independent components:\\
	\begin{minipage}{0.45\columnwidth}
		\vspace{-0.3cm}
		\begin{center}
	\begin{eqnarray*}
		\pmat{E} &= \sum_{i=1}^{Q_{E}} \theta_{i}^{E}(\mu) \textbf{E}_{i},\\
		\pmat{B} &= \sum_{i=1}^{Q_{B}} \theta_{i}^{B}(\mu) \textbf{B}_{i},
	\label{eq:aff1}
	\end{eqnarray*}
		\end{center}
	\end{minipage}
	\hspace{0.02\columnwidth}
	\begin{minipage}{0.45\columnwidth}
	\vspace{-0.3cm}
		\begin{center}
	\begin{eqnarray*}
		\pmat{A} &= \sum_{i=1}^{Q_{A}} \theta_{i}^{A}(\mu) \textbf{A}_{i},\\
		\pmat{C} &= \sum_{i=1}^{Q_{C}} \theta_{i}^{C}(\mu) \textbf{C}_{i}.
	\label{eq:aff2}
	\end{eqnarray*}
		\end{center}
	\end{minipage}
	\newline \\
In the offline stage, the expensive parameter independent terms are evaluated once and stored for repeated use in the online stage. 

The standard process of the RB method for time-dependent parametric systems is presented in~\cref{alg:1} which is known as the POD-Greedy algorithm. 
It constructs the projection matrix $\textbf{V}$ by selecting a set of samples $\{ \mu_{i} \}_{i=1}^{r}$ from a given parameter training set, through a greedy algorithm. The FOM solutions $x^{k}(\mu)$, $k = 1,\ldots,K$, at a newly selected sample $\mu_{i}$ are assembled into a matrix $\textbf{X} \coloneqq [x^{1}(\mu_{i}), \ldots, x^{K}(\mu_{i}) ]$. The matrix $\textbf{V}$ is iteratively enriched by projecting the solution vectors in $\mathbf{X}$ onto the current subspace $\text{range}(\mathbf{V})$, see Steps 3 - 5 in~\cref{alg:1}. Crucial to the success of the RB method is the \textit{a posteriori} error estimator,
	\begin{align*}
		\Delta^{k+1}(\mu) \geq \| y^{k+1}(\mu) - y_{r}^{k+1}(\mu) \|.
	\end{align*}
It quantifies the error between the FOM and ROM for each parameter, without having to solve the FOM.
The parameters at which the FOM is to be solved to enrich $\mathbf{V}$, is picked iteratively by maximizing the error estimator over the training set. Although \cref{eq:rom} is of reduced order, the complexity of evaluating the nonlinear term $f_{r}(\textbf{V} x_{r}^{k}(\mu); \mu)$ remains of order $\mathscr{O}(n)$. This is the so-called \emph{lifting bottleneck}. Techniques like Empirical Interpolation Method (EIM) or its discrete variation, the Discrete Empirical Interpolation Method (DEIM), can be used to obtain a ROM whose computational complexity is independent of the full order $n$. For more details on this approach, we refer the reader to \cite{morBarMNetal04,morChaS10}. 
In the following section we discuss the idea of adaptive parameter sampling.
\section{Adaptive Parameter Sampling}
\label{sec:3}
The quality of the ROM depends on how well the training set represents the parameter domain. In the standard approach, if a finely sampled training set is used for a parameter domain with high-dimension, it can be computationally expensive. If otherwise, a coarse training set is used, it will not guarantee that a reliable ROM can be derived. Consequently, one faces the situation of trial and error. The approach we propose here is based on a surrogate error model generated through RBF interpolation. It aims to overcome the drawbacks of the standard approach using a fixed training set.
\subsection{Radial Basis Functions}
\label{subsec:rbf}
RBF interpolation is an efficient way to approximate scattered data in high-dimensional space. Consider the domain $\Omega \subseteq \mathbb{R}^{d}$ and the function $h \coloneqq \mathbb{R}^{d} \xrightarrow{} \mathbb{R}$. Assume that $h$ is difficult to evaluate, or known only at a few points in $\Omega$. Let $\Phi$ be the \emph{kernel function} having the property that $\Phi(\mu_{1}, \mu_{2}) = \Phi(\| \mu_{1} - \mu_{2} \|), \forall \mu_{1}, \mu_{2} \in \Omega$. Typically, the Euclidean norm is used. Such kernels have the name of \emph{radial basis functions} since they are radially symmetric. We seek an approximation $s: \mathbb{R}^{d} \xrightarrow{} \mathbb{R}$ to the function $h$ given as
	\begin{equation}
		s(\mu) := \sum_{i=1}^{\ell} c_{i} \Phi(\| \mu - \mu_{i} \|), \quad \forall \mu\in \Omega.
	\label{eq:rbf}
	\end{equation}
The coefficients $\{ c_{i} \}_{i=1}^{\ell}$ are determined by enforcing the interpolation condition $h(\mu_{i}) = s(\mu_{i}), \, i = 1, 2, \ldots, \ell$. This amounts to solving the linear system of equations
	\begin{equation}
		\underbrace{\begin{bmatrix}
			\Phi(\mu_{1}, \mu_{1}) & \Phi(\mu_{1}, \mu_{2}) & \cdots & \Phi(\mu_{1}, \mu_{\ell}) \\
			\Phi(\mu_{2}, \mu_{1}) & \Phi(\mu_{2}, \mu_{2}) & \cdots & \Phi(\mu_{2}, \mu_{\ell}) \\
			\vdots & \vdots & \ddots & \vdots \\
			\Phi(\mu_{\ell}, \mu_{1}) & \Phi(\mu_{\ell}, \mu_{2}) & \cdots & \Phi(\mu_{\ell}, \mu_{\ell})
		\end{bmatrix}}_{\mathbf{R}}
		\underbrace{\begin{bmatrix}
			c_{1} \\ c_{2} \\ \vdots \\ c_{\ell}
		\end{bmatrix}}_{c}
		=
		\underbrace{\begin{bmatrix}
			h(\mu_{1})\\ h(\mu_{2}) \\ \vdots \\ h(\mu_{\ell})
		\end{bmatrix}}_{\mathfrak{d}}
		\label{eq:rbfsystem}
	\end{equation}
	
The points $\{\mu_{i}\}_{i=1}^{\ell} \in \Omega$ are called the \emph{centers} of the radial basis functions. Under the assumption that the centers are pairwise distinct, it can be proven that the RBF kernel matrix $\mathbf{R}$ is positive definite for some suitable choice of radial basis functions $\Phi$ and thus \cref{eq:rbfsystem} has a unique solution. We refer to \cite[Chapter 6]{Wed05} for the proof and detailed treatment. \cref{tab:kernels} provides a list of commonly used radial basis functions. All of the radial basis functions in \cref{tab:kernels} have global support. There exist also radial basis functions with local support. For more details we refer the reader to \cite{Buh03}.

In practice, the class of radial basis functions that have the positive definite property ($\mathbf{R}$ being positive definite) is limited to a few, such as Gaussian and Inverse multiquadric. The class of admissible basis functions ($\Phi(\cdot)$) can be expanded by defining the so called \emph{conditionally positive definite} functions, by which the positive definiteness can be satisfied by imposing some additional constraints given as,
\begin{equation*}
\sum_{j=1}^{M} c_{i} p_{j}(\mu) = 0, \quad i = 1, 2, \ldots , \ell.
\end{equation*}
The functions $p_{1}, p_{2}, \ldots, p_{M}$  are a basis of the polynomial space with suitable degree. In practice, we choose $M$ to be equal to the number of scalar parameters $d$. With the new conditions imposed, the radial basis interpolant now becomes,
	\begin{equation}
		s(\mu) := \sum_{i=1}^{\ell} c_{i} \Phi(\| \mu - \mu_{i} \|) + \sum_{j=1}^{M} \lambda_{j} p_{j}(\mu).
	\end{equation}
	\begin{table}[t!]
		\scriptsize
		\centering
		\begin{tabular}{|c|c|}
			\hline
			Kernel\Tstrut & $\Phi(r) := \Phi(\| \mu - \mu_{i} \|)$\Tstrut \\[0.5ex]
			\hline
			\rule{0pt}{3ex} Gaussian & $e^{- \sigma^{2} r^{2}}$ \\[0.5ex] 
			Thin-plate splines & $r^{2} \ln(r)$\\[1ex]	 
			Multiquadric & $\sqrt{r^{2} + \sigma^{2}}$ \\[1ex] 
			Inverse multiquadric & $\dfrac{1}{\sqrt{r^{2}+ \sigma^{2}}}$\\[2ex] \hline
		\end{tabular}
	\caption{Common radial basis kernels.}
	\label{tab:kernels}
	\end{table}
We then obtain a saddle-point system of dimension $N_{\text{RBF}} := (M + \ell) \times (M + \ell)$:
	\begin{equation}
		\begin{bmatrix}
			\mathbf{R} & \mathbf{P}\\
			\mathbf{P}^{T} & 0
		\end{bmatrix}
		\begin{bmatrix}
			c \\ \lambda
		\end{bmatrix}
		=
		\begin{bmatrix}
			\mathfrak{d} \\ 0
		\end{bmatrix}
	\label{eq:rbfsystem-saddle}
	\end{equation}
With a proper choice of $p_{1}, p_{2}, \ldots, p_{M}$, the augmented coefficient matrix is positive definite for all choices of radial basis functions in \cref{tab:kernels} and this ensures the uniqueness of the interpolant. For a detailed discussion of the rationale behind the idea of conditionally positive definite functions, we refer to \cite[Chapter 6]{Wed05}. 
\subsection{POD-Greedy algorithm with adaptive sampling}
\label{subsec:POD_greedy_adapt}
Now we detail the proposed approach of adaptively constructing the training set for the RB method. 
The proposed adaptive parameter sampling scheme is general and can be combined with the standard greedy algorithm for steady parametric systems, or with the POD-greedy algorithm (\cref{alg:1}) for time-dependent systems, to adaptively update the training set and reduce the offline costs.  Since the framework is similar, we only present the adaptive parameter sampling approach with POD-Greedy for time-dependent problems. The standard POD-Greedy algorithm is presented in~\cref{alg:1}. POD-greedy with adaptive sampling is presented as~\cref{alg:1_ad}, where it can be seen that the training set ($\Xi_c$) is updated at each iteration. 

\cref{alg:1_ad} is initialized using a coarse training set $\Xi_{c}$, with cardinality $N_{c}$. At the end of each iteration, we compute the error estimator {\it only} over the {\it coarse} training set $\Xi_c$, i.e., $\Delta(\mu), \, \forall \mu \in \Xi_{c}$ (compare Step 6 in~\cref{alg:1} with Step 4 in~\cref{alg:1_ad}). We use this as the \emph{input data} to form the radial basis interpolant over the fine training set, i.e., $s(\mu)\,, \, \forall \mu \in \Xi_{f}$, where $\Xi_{f}$ is a finely sampled training set of cardinality $N_{f} \gg N_{c}$. The evaluation costs are $\mathcal{O}((N_{c}+M)^{3})$ for identifying the interpolant coefficients, and $\mathcal{O}(N_{f}N_{c})$ for interpolating over the fine training set. Here $M$ is the number of  polynomial functions added to make the augmented coefficient matrix in \cref{eq:rbfsystem-saddle} conditionally positive definite. Since $N_{c}, M$ are both small numbers, the cost of evaluating  $s(\mu)$ is much smaller than the cost of computing $\Delta(\mu)$ which includes solving the ROM for each $\mu$. We enrich the coarse training set by adding new parameters from $\Xi_{f}$ to $\Xi_{c}$. We add those $n_{\text{add}}$ new parameters that have the largest magnitude of $s(\mu)$, where $n_{\text{add}}$ can be fixed and user-defined. One heuristic way of varying $n_{\text{add}}$ adaptively is $n_{\text{add}} = \log_{10} \bigg( \Biggl\lfloor \dfrac{\max \limits_{\mu \in \Xi_{c}} \Delta(\mu)} {\texttt{tolerance}} \Biggr\rfloor \bigg).$
Additionally, we also monitor the coarse training set to identify those points $\breve{\mu} \in \Xi_{c}$ with $\Delta(\breve{\mu}) < \texttt{tol}$. Those points are then removed from $\Xi_{c}$, meaning that their corresponding full solutions need not be computed to enrich the RB basis. In this way, the coarse training set gets updated iteratively and its size remains as small as possible, avoiding many unnecessary full simulations of the FOM. 

In this work, we make use of the estimator proposed in \cite{morCheFB19a}. It is given by,
	\begin{equation*}
		\Delta^{k+1}(\mu) := \bar{\Psi} \| r_{\text{pr}}^{k+1} \|,
	\end{equation*}
where $r_{\text{pr}}^{k+1}$ is the residual obtained by substituting the approximate solution $\mathbf{V} \mathbf{x}_{r}^{k+1}$ in \cref{eq:fom}, at each time instance $t^{k+1}$. The term $\bar{\Psi}$ involves the norm of the approximate solution to the dual system and its residual. Additionally, it involves two terms that need to be estimated. The first is the \textit{inf-sup} constant, given by the smallest singular value of the matrix $\mathbf{E}(\mu)$ (when considering the Euclidean norm). The second term to be estimated involves the residual of an auxiliary system. For a more detailed discussion, we refer the reader to \cite{morCheFB19a}. In the implementation of Algorithms 3 and 4, we have made use of the method proposed in \cite{morManN15} to estimate the inf-sup constant efficiently.

\renewcommand{\algorithmicrequire}{\textbf{Input:}}
\renewcommand{\algorithmicensure}{\textbf{Output:}}
\begin{algorithm}[t]
	\caption{Standard POD-Greedy}
	\begin{algorithmic}[1]
		\label{alg:1}
		\REQUIRE Training set $\Xi$, Tolerance (\texttt{tol}).
		
		\ENSURE $\mathbf{V}$
		
		\STATE Initialize. $\mathbf{V} = [\,]$, $\mu^{*}\in \Xi$. 
		
		\WHILE{$\Delta(\mu^{*})$ $>$ \texttt{tol}} 
		
		\STATE Compute full order solution at $\mu^{*}$. $\mathbf{X} = [x(t^{1}, \mu^{*}), \, x(t^{2}, \mu^{*}), \, \ldots, \, x(t^{K}, \mu^{*})]$.
		
		\STATE Form $\bar{\mathbf{X}} := \mathbf{X} - \mathbf{V} \mathbf{V}^{T} \mathbf{X}$; Obtain $\bar{\mathbf{X}} \xrightarrow{\texttt{svd}}$ $\mathbf{U} \Sigma \mathbf{W}^{T}$. 
		
		\STATE Update $\mathbf{V}$:= $\texttt{orth}[ \mathbf{V}$, $\mathbf{U}(: \, , 1)]$.		
		
		\STATE Solve the ROM and compute the error estimator $\Delta(\mu)$, $\forall \mu \in \Xi$.
		
		\STATE Next $\mu^{*}$ is chosen as,
		$\mu^{*}$ := arg $\max\limits_{\mu \in \Xi} \Delta(\mu)$.		
		\ENDWHILE
	\end{algorithmic}
\end{algorithm}	
	\renewcommand{\algorithmicrequire}{\textbf{Input:}}
		\renewcommand{\algorithmicensure}{\textbf{Output:}}
		\begin{algorithm}[h]
		\caption{POD-Greedy with adaptive parameter sampling}
		\begin{algorithmic}[1]
		\label{alg:1_ad}
		\REQUIRE Coarse training set $\Xi_{c}$, Fine training set $\Xi_{f}$, Tolerance (\texttt{tol}).
		
		\ENSURE $\mathbf{V}$.
		
		\STATE Initialize. $\mathbf{V} = [\,]$, $\mu^{*} \in \Xi_{c}$. 
		
		\WHILE{$\Delta(\mu^{*})$ $>$ \texttt{tol}} 
		
			\STATE Perform Steps 3-5 in~\cref{alg:1}.
			
			\STATE Solve the ROM and compute the error estimator $\Delta(\mu)$, $\forall \mu \in \Xi_{c}$. 

			\STATE Form RBF interpolant $s(\mu)$ of $\Delta(\mu)$ over $\Xi_{f}$.
			
			\STATE Calculate $n_{\text{add}}$, pick $\{\mu_{1}, \ldots, \mu_{\text{add}} \}$ from $\Xi_{f}$ with largest errors measured by $s(\mu)$.
			\STATE Update the training set. $\Xi_{c} \coloneqq [\Xi_{c} \cup \{\mu_{1}, \ldots, \mu_{\text{add}} \}]$ (Add samples to the current training set $\Xi_c$). Find all $\breve{\mu}_{i}, \forall i = 1,\ldots, \text{n}_{\text{del}}$ with $\Delta(\breve{\mu}_{i}) < \texttt{tol}$. Set $\Xi_{c} \coloneqq \Xi_{c} \setminus \{ \breve{\mu}_{i}, \forall i = 1,\ldots, \text{n}_{\text{del}} \}$ (Remove unnecessary samples from the current training set $\Xi_c$).
			\STATE Next $\mu^{*}$ is chosen as,
			$\mu^{*}$ := arg $\max\limits_{\mu \in \Xi_{c}} \Delta(\mu)$.
					
		\ENDWHILE
		\end{algorithmic}
		\end{algorithm}	 
\begin{remark}[Interpolating small values]
					The magnitudes of the errors we interpolate are often very small, especially at the latter iterations of the greedy algorithm. To ensure good, stable interpolation, we consider the input data to the RBF as the base $10$ logarithm of the estimated errors. After interpolation, we perform the anti-logarithm and project the logarithm value back to the actual error value.
				\end{remark}
\subsection{A fully adaptive POD-Greedy-(D)EIM algorithm}
\label{subsec:fulladapt}
An adaptive version of~\cref{alg:1}, called the Adaptive POD-Greedy-(D)EIM algorithm, is proposed in \cite{morCheFB19a} for nonlinear systems, which aims to update the RB, (D)EIM basis with an adaptively chosen number of basis vectors at each iteration of the greedy algorithm. In case of a nonlinear or nonaffine system, the standard POD-Greedy algorithm usually precomputes the (D)EIM basis and the interpolation index matrix, $\mathbf{U}_{f}$ and $\mathbf{P}_{f}$. Unlike the standard POD-Greedy algorithm (\cref{alg:1}), the Adaptive POD-Greedy algorithm (\cref{alg:2}) updates the (D)EIM matrix inside the greedy algorithm. This avoids many additional FOM simulations caused by the separate (D)EIM computation outside the POD-Greedy loop. To ensure a reliable ROM, an efficient a posteriori error estimator ($\Delta(\mu)$) is used.

The algorithm updates the RB vectors and (D)EIM basis vectors at every iteration (Steps 5-8 in~\cref{alg:2}).
The number of RB basis vectors, $\delta_{\text{RB}}$, and that of (D)EIM basis vectors, $\delta_{\text{DEIM}}$, to be added to/removed from the current basis are computed adaptively, see Step 8 in~\cref{alg:2}. For details of computing $\delta_{\text{RB}}$ and $\delta_{\text{DEIM}}$, we refer the reader to~\cite{morCheFB19a}. We note that for this algorithm, the training set is fixed.

\renewcommand{\algorithmicrequire}{\textbf{Input:}}
\renewcommand{\algorithmicensure}{\textbf{Output:}}
\begin{algorithm}[t!]
	\caption{Adaptive POD-Greedy-(D)EIM algorithm}
	\begin{algorithmic}[1]
		\label{alg:2}
		\REQUIRE Training set $\Xi$, Tolerance (\texttt{tol}).
		
		\ENSURE $\mathbf{V}$, DEIM matrices $\mathbf{U}_{f}, \mathbf{P}_{f}$.
		
		\STATE Initialize. $\mathbf{V} = [\,]$, $\mu^{*} \in \Xi$, $\delta_{\text{RB}} = 1$, $\ell_{\text{DEIM}} = 1$, $\delta_{\text{DEIM}}=0$, $\mathbf{U}_{f} = [\,], \mathbf{P}_{f} = [\,]$, $\mathbf{F} = [\,]$. 
		
		\WHILE{$\Delta(\mu^{*})$ $>$ \texttt{tol}} 
		
			\STATE Compute full order solution at $\mu^{*}$. Form the snapshot matrices: \\ $\mathbf{X} = [x(t^{1}, \mu^{*}), \, x(t^{2}, \mu^{*}), \, \ldots, \, x(t^{K}, \mu^{*})]$ and \\ $\mathbf{F}^{\text{new}} = [f(x(t^{1}, \mu^{*})), \, f(x(t^{2}, \mu^{*})), \, \ldots, \, f(x(t^{K}, \mu^{*}))]$.
			
			\STATE Form $\bar{\mathbf{X}} := \mathbf{X} - \mathbf{V} \mathbf{V}^{T} \mathbf{X}$; obtain $\bar{\mathbf{X}} \xrightarrow{\texttt{svd}}$ $\mathbf{U} \Sigma \mathbf{W}^{T}$. 
			\STATE Update $\mathbf{V}$:= $\texttt{orth}\big([ \mathbf{V}$, $\mathbf{U}(: \, , 1:\delta_{\text{RB}})]\big)$ if $\delta_{\text{RB}} \geq0$; remove $\delta_{\text{RB}}$ vectors from 
	 $\mathbf{V}$ if $\delta_{\text{RB}} <0$. 
			
			\STATE Compute $\ell_{\text{DEIM}}$ (D)EIM interpolation basis vectors from snapshots of the nonlinear function, $\mathbf{F}: = [\mathbf{F}, \mathbf{F}^{\text{new}}]$.
			
			\STATE Solve ROM and compute the error estimator $\Delta(\mu)$, $\forall \mu \in \Xi$.

			\STATE Based on $\Delta(\mu)$, decide the new number $\delta_{\text{RB}}$ and $\delta_{\text{DEIM}}$ using the adaptive scheme in~\cite{morCheFB19a}. $\ell_{\text{DEIM}}=\ell_{\text{DEIM}}+\delta_{\text{DEIM}}$.
			
			\STATE Next $\mu^{*}$ is chosen as,
			$\mu^{*}$ := arg $\max\limits_{\mu \in \Xi} \Delta(\mu)$.		
		\ENDWHILE
	\end{algorithmic}
\end{algorithm}	
With the proposed adaptive sampling technique, we present a fully adaptive POD-Greedy-(D)EIM algorithm: \cref{alg:3}. It is a combination of~\cref{alg:2} with the adaptive sampling approach as presented in~\cref{alg:1_ad}. Since the adaptive sampling scheme is the same as in~\cref{alg:1_ad}, we omit the detailed explanation for~\cref{alg:3}. One only needs to keep in mind that the error estimator is computed {\it only} over the {\it coarse} training set, see Step 4. To estimate the errors at samples in the fine training set, the error surrogate $s(\mu)$ is computed and used, see Step 6. 
		\renewcommand{\algorithmicrequire}{\textbf{Input:}}
		\renewcommand{\algorithmicensure}{\textbf{Output:}}
		\begin{algorithm}[t!]
		\caption{Adaptive POD-Greedy-(D)EIM with adaptive parameter sampling}
		\begin{algorithmic}[1]
		\label{alg:3}
		\REQUIRE Coarse training set $\Xi_{c}$, Fine training set $\Xi_{f}$, Tolerance (\texttt{tol}).
		
		\ENSURE $\mathbf{V}$, DEIM basis $\mathbf{U}_{f}, \mathbf{P}_{f}$.
		
		\STATE Initialize. $\mathbf{V} = [\,]$, $\mu^{*}\in \Xi_c$, $\ell_{\text{RB}} = 1$, $\ell_{\text{DEIM}} = 1$, $\delta_{\text{DEIM}} = 0$, $\mathbf{U}_{f} = [\,], \mathbf{P}_{f} = [\,]$. 
		
		\WHILE{$\Delta(\mu^{*})$ $>$ \texttt{tol}} 
		
			\STATE Perform Steps 3-6 in~\cref{alg:2}.

			\STATE Solve ROM and compute the error estimator $\Delta(\mu)$, $\forall \mu \in \Xi_{c}$ . 
			\STATE Based on $\Delta(\mu)$, decide on the new number $\delta_{\text{RB}}$ and $\delta_{\text{DEIM}}$ using the adaptive scheme in~\cite{morCheFB19a}. $\ell_{\text{DEIM}}=\ell_{\text{DEIM}}+\delta_{\text{DEIM}}$.

			\STATE Perform Steps 5-8 in~\cref{alg:1_ad}

		\ENDWHILE
		\end{algorithmic}
		\end{algorithm}	
\subsection{A strategy for shape parameter selection}
\label{subsec:RBF_shape}

In \cite{morTaineA15}, the authors propose a technique to adaptively sample parameters based on Kriging interpolation, where several hyper-parameters need to be estimated. As an alternative, we use RBF to interpolate an efficient output error estimator. Based on the choice of the kernel used, RBF involves at most one free parameter the user has to specify. This is the case for Gaussian or Multiquadric kernels, where the shape parameter $\sigma$ needs to be determined. For the case of Thin-plate spline kernels, there are no free parameters to choose.
A good choice of the shape parameter is essential to ensure that the RBF kernel matrix ($\mathbf{R}$) remains well-conditioned~\cite{FasZ07}.
In~\cite{Rip99}, a heuristic approach is introduced to determine the optimal shape parameter based on the idea of the Leave One Out Cross-Validation (LOOCV) strategy commonly used in the field of statistics~\cite{Jametal13}. The main idea of LOOCV is to utilize the available data and find a $\sigma$ that best fits the data. Below, we briefly describe how we use it to estimate the shape parameter. For more details on the method and its generalization - the k-fold cross validation, we refer to \cite{Jametal13}.
	\begin{enumerate}[(i)]
		\item We start with the available $\ell$ centers and data points $\big( \mu_{i}, \Delta(\mu_{i}) \big), \forall i = 1, \ldots, \ell$ for the RBF interpolation.
		\item For every $\mu_i, \forall i = 1, \ldots, \ell$, compute a `less accurate' radial basis interpolant $s^{\mu_i}(\mu)$ by removing the ${i}^{\text{th}}$ row and ${i}^{\text{th}}$ column of the RBF kernel matrix $\mathbf{R}$ and the ${i}^{\text{th}}$ row of the input data vector $\mathfrak{d}$. Note that $s^{\mu_i}(\mu_i)$ is actually an approximation of $h(\mu_i)$, whereas,
$s(\mu_i)=h(\mu_i)$ (recall~\cref{eq:rbfsystem}).
 		\item Following this, we have the error between $h(\mu_i)$ and $s^{\mu_i}(\mu_i)$: $e^{\mu_{i}}:=h(\mu_i)-s^{\mu_i}(\mu_i)$. Since $h(\mu_i)=s(\mu_i)$, we have
			\begin{align*}
				e^{\mu_{i}} = s(\mu_{i}) - s^{\mu_{i}}(\mu_{i}), \quad  i = 1, \ldots, \ell.
			\end{align*}
		\item Form the error vector $\mathbf{e} = [e^{\mu_{1}}, e^{\mu_{2}}, \ldots, e^{\mu_{\ell}}]$.
		\item Choose the optimal $\sigma$ as $\sigma^{*} = \arg \min\limits_{\sigma \in [\sigma_{\text{min}}\,,\sigma_{\text{max}}]} \| \mathbf{e}(\sigma) \|_{2}$.
	\end{enumerate}
We solve the minimization problem using the MATLAB\textsuperscript \textregistered\,function \texttt{fminbnd}, as suggested in \cite{FasZ07}. For the success of the adaptive sampling approach, the following criteria are crucial:
	\begin{enumerate}
		\item \emph{Good data:} The input data to the RBF should be `good'. This means that we are interpolating something meaningful to the problem we are trying to approximate.
		\item \emph{Good shape parameter:} The RBF interpolation should be robust, ensuring that the interpolated values can act as `good' surrogates of the actual values.
	\end{enumerate}
We ensure the first criterion by using an efficient output error estimator for the RB method from \cite{morCheFB19a}. This is a residual based estimator suitable for nonlinear dynamical systems and tailored for output error estimation. 
Meanwhile, the second criterion is ensured by adopting the LOOCV strategy described above.
\section{Numerical results}
\label{sec:4}
In this section, we validate the proposed adaptive parameter sampling method using three examples. The first is a nonlinear Burgers' equation model with one parameter. The second is a two-parameter linear convection-diffusion model. The last one is a three-parameter model of a microthruster unit. All numerical tests were performed in MATLAB\textsuperscript \textregistered 2015a, on a laptop with Intel\textsuperscript \textregistered Core\textsuperscript \texttrademark i5-7200U @ 2.5 GHZ, with 8 GB of RAM. 
For all the examples,
	\begin{itemize}

	\item The Adaptive POD-Greedy-(D)EIM algorithm proposed in \cite{morCheFB19a} is used for basis construction. Therefore, we compare the results of~\cref{alg:2} with fixed training set with those of~\cref{alg:3} with adaptive training set. 

	\item The maximal output error over the test set $\Xi_{\text{test}}$ is defined as
		\begin{equation}
			\begin{aligned}
				\varepsilon_{\text{max}}:=\max\limits_{\mu \in \Xi_{\text{test}}} \bigg( \frac{1}{K}\sum_{k=1}^{K} \| y^{k}(\mu) - \bar{y}^{k}(\mu) \| \bigg).
			\end{aligned}
			\label{eq:errorplot}
		\end{equation}

	\item The mean (over time steps) output error over $\Xi_{\text{test}}$ is defined as
		\begin{equation}
			\begin{aligned}
			\varepsilon(\mu):=\bigg( \frac{1}{K}\sum_{k=1}^{K} \| y^{k}(\mu_i) - \bar{y}^{k}(\mu_i) \| \bigg), \quad \forall \mu \in \Xi_{\text{test}}.
			\end{aligned}
			\label{eq:testerrorplotfinal}
		\end{equation}

	\item The semi-implicit Euler scheme is used for time discretization.

	\item If not particularly pointed out, the {\it initial} coarse training set $\Xi_c$ for~\cref{alg:3} is the same as the fixed training set $\Xi$ for~\cref{alg:2}.

	\end{itemize}
\subsection{Burgers' equation} 
			\newlength\fheight
			\newlength\fwidth
We consider the $1$-D model of the viscous Burgers' equation. The PDE is defined in the spatial domain $w \in \Omega \coloneqq [0\,,1]$ and for time $T \coloneqq [0\,, 2]$ as
	\begin{equation*}
\begin{array}{l}
		\frac{dv}{dt} + v \frac{\partial v}{\partial w}= q \frac{\partial^{2} v}{\partial w^{2}} + s(w,t),\quad
v(0, t) = 0,\quad \frac{\partial v(1, t)}{\partial w} = 0.
	\end{array}
\end{equation*}
where $v(w,t)$ is the unknown variable and $q \in \mathscr{P} := [0.001\,,1]$ is the viscosity parameter. 
The initial condition is set as zero and a constant input to the system is  $s(w,t) \equiv 1$. The output is monitored on the last spatial point in the domain $y = v(1,t)$. The model has $N = 500$ equations after discretization in space. 

We show the results of~\cref{alg:2} and~\cref{alg:3}, respectively. The tolerance is set as $\texttt{tol} = 10^{-5}$. To implement the RB method without adaptive parameter sampling, i.e.~\cref{alg:2}, we choose an initial training set consisting of 10 random parameter samples, using the \texttt{rng} command in MATLAB\textsuperscript \textregistered, with the seed value $112$ for the random number generator $\texttt{twister}$.
The fine training set $\Xi_{f}$ for~\cref{alg:3} consists of 300 random samples from the same parameter domain, sampled using the seed value at $114$ for the random number generator \texttt{twister}. Finally, 100 different random samples are considered for the test parameter set $\Xi_{\text{test}}$, sampled using \texttt{simdTwister} with the seed value $200$. To construct the error surrogate model $s(\mu)$ in~\cref{alg:3}, the Inverse Multiquadric (IMQ) kernel function is used for the RBF interpolation. Cross validation LOOCV has been applied to specify the shape parameter.
Since the model is nonlinear, DEIM \cite{morChaS10} is used in order to efficiently compute the nonlinear term. 

In~\cref{fig:1a} shows the decay of the maximal output error, i.e. $\varepsilon_{\text{max}}$, over the test set $\Xi_{\text{test}}$. It is
calculated at every iteration of either~\cref{alg:2} or~\cref{alg:3}. The former results in a ROM of dimension $11$, while the latter results in a ROM of dimension $28$. Evidently, using an adaptive training set leads to a better convergence. \cref{alg:2} converges in 7 iterations, while~\cref{alg:3} requires 18 iterations to converge. However, as seen from the output error $\varepsilon(\mu)$ of the final ROM over
the test set $\Xi_{\text{test}}$ in~\cref{fig:1b},~\cref{alg:3} produces a ROM with errors being uniformly below the tolerance, whereas, the ROM computed using~\cref{alg:2} has large errors above the tolerance. In \cref{fig:2}, we show the evolution of the coarse training set $\Xi_{c}$ for different iterations of \cref{alg:3}. It is seen that the algorithm tends to pick parameters in the low-viscosity regions (close to 0.001), as expected, since the solution of the PDE tends to be `less smooth', requiring more basis functions to approximate.	
			\setlength\fheight{3cm}
			\setlength\fwidth{3cm}		
			\begin{figure}[t]
				\centering
				\subfloat[Convergence of~\cref{alg:2} (fixed) and~\cref{alg:3} (adaptive).]{\label{fig:1a}
%
%
\begin{tikzpicture}

\begin{axis}[%
width=\fwidth,
height=\fheight,
at={(0\fwidth,0\fheight)},
scale only axis,
xmin=1,
xmax=18,
xlabel={$i$th iteration},
xlabel style={font=\color{white!15!black},font = \small},
ymode=log,
ymin=1e-07,
ymax=1,
yminorticks=true,
ylabel={$\varepsilon_{\text{max}}$},
ylabel style={font=\color{white!15!black}, font = \small},
legend style={legend cell align=left, align=left, draw=white!15!black, font = \scriptsize},
axis background/.style={fill=white}
]
\addplot [color=mpiblue, very thick]
  table[row sep=crcr]{%
1	0.0749162805766902\\
2	0.00126429197278156\\
3	0.000612268764195509\\
4	0.00023741747749495\\
5	0.000131460188830449\\
6	0.000194452625010194\\
7	9.45836294893178e-05\\
};
\addlegendentry{Fixed}
\addplot [color=mpired, thick, densely dashed]
  table[row sep=crcr]{%
1	0.0749162805766902\\
2	0.00198115745569395\\
3	0.00021486974178479\\
4	7.56787710711788e-05\\
5	8.32904628397757e-05\\
6	1.8555378641312e-05\\
7	3.01739948583568e-05\\
8	9.14040147564984e-05\\
9	0.000117369205152922\\
10	2.31402970384658e-05\\
11	2.26584022445925e-06\\
12	1.74878041829339e-06\\
13	2.26966081879933e-06\\
14	3.64109365016241e-06\\
15	8.64198585901602e-07\\
16	2.40274783883083e-07\\
17	1.1639542267462e-07\\
18	1.31255334646564e-07\\
};
\addlegendentry{Adaptive}
\end{axis}
\end{tikzpicture}%
} \hfill
				\subfloat[Error $\varepsilon(\mu)$ of the final ROM over $\Xi_{\text{test}}$.]{\label{fig:1b}
%
%
\begin{tikzpicture}

\begin{axis}[%
width=\fwidth,
height=\fheight,
at={(0\fwidth,0\fheight)},
scale only axis,
xmin=0,
xmax=100,
xlabel={$i$th parameter in $\Xi_{\text{test}}$},
xlabel style={font=\color{white!15!black},font = \small},
ymode=log,
ymin=1e-09,
ymax=0.001,
yminorticks=true,
ylabel={$\varepsilon(\mu_i)$},
ylabel style={font=\color{white!15!black}, font = \small},
legend style={legend cell align=left, align=left, draw=white!15!black, font = \scriptsize},
axis background/.style={fill=white}
]
\addplot [color=mpiblue, thick]
  table[row sep=crcr]{%
1	1.63404077285713e-07\\
2	2.42328269769726e-07\\
3	1.5129111081653e-07\\
4	2.08544637266294e-07\\
5	2.36712946954469e-07\\
6	2.43866484698251e-07\\
7	1.96846793179659e-07\\
8	1.64959853447615e-07\\
9	2.17519884063291e-07\\
10	2.30248761151276e-07\\
11	2.20394656639462e-07\\
12	1.10357863011848e-06\\
13	1.90498768557631e-07\\
14	1.98606752660235e-07\\
15	2.26261300104468e-07\\
16	2.09195734067088e-07\\
17	9.04866655656654e-08\\
18	1.25273996305494e-06\\
19	3.0932095859611e-07\\
20	2.25762400887595e-07\\
21	2.33949916830326e-07\\
22	1.82235523303053e-07\\
23	7.91864804617065e-07\\
24	9.09353362644346e-08\\
25	2.4217081111277e-05\\
26	1.54239364295436e-07\\
27	1.66941353060562e-07\\
28	1.92301207948319e-07\\
29	1.18005910521772e-06\\
30	2.42253462416739e-07\\
31	2.41166235214136e-07\\
32	1.6986869608051e-07\\
33	2.18227768025251e-07\\
34	2.44335693901605e-07\\
35	5.37739205921492e-07\\
36	9.89601304317055e-07\\
37	1.41642279502802e-07\\
38	2.24692734701223e-07\\
39	1.09888233852322e-06\\
40	1.56978879231144e-07\\
41	7.02298452070734e-05\\
42	2.32329957522562e-07\\
43	1.15860386123064e-07\\
44	9.35349645367644e-07\\
45	1.19626454992174e-06\\
46	9.45836294893178e-05\\
47	1.77916848414081e-07\\
48	2.86860956176255e-07\\
49	7.7883185722964e-07\\
50	2.28093901337589e-05\\
51	1.93922748373946e-07\\
52	2.37085724973026e-07\\
53	2.33931526363073e-07\\
54	2.44709535069186e-07\\
55	1.48431873724478e-07\\
56	1.75619976253884e-07\\
57	1.04850782188526e-07\\
58	8.4681403641196e-07\\
59	1.75571472036902e-07\\
60	8.33853080536823e-07\\
61	2.27009384245126e-07\\
62	1.86049311565102e-07\\
63	2.27974465200908e-07\\
64	2.19665411226088e-07\\
65	1.96317215656485e-07\\
66	1.23373540531276e-06\\
67	2.31277220172518e-07\\
68	2.40678101077539e-07\\
69	1.54420457168568e-07\\
70	2.4416723173337e-07\\
71	2.2721975310553e-07\\
72	1.36652165718016e-05\\
73	9.93951837096945e-08\\
74	1.50688144237903e-07\\
75	1.65622442944228e-07\\
76	2.84754217769299e-07\\
77	2.41847078531727e-07\\
78	2.32190544212417e-07\\
79	1.97686700537375e-07\\
80	8.75345505877575e-07\\
81	1.59596524470599e-07\\
82	2.1755189172207e-07\\
83	1.24794038674367e-06\\
84	1.64797395887411e-07\\
85	3.11374478524268e-05\\
86	6.60983096234728e-06\\
87	2.27530286816067e-07\\
88	2.43902110421483e-07\\
89	7.09162819340408e-07\\
90	2.27945959312491e-07\\
91	2.44190906349559e-07\\
92	1.74061138678795e-07\\
93	9.24144856329188e-08\\
94	1.47278617542886e-07\\
95	2.35063388704806e-07\\
96	2.20537928678711e-07\\
97	9.06496554323667e-07\\
98	1.4579772020621e-07\\
99	2.36695420634885e-07\\
100	1.05758219700504e-07\\
};
\addlegendentry{Fixed}
\addplot [color=mpired, thick,densely dashed]
  table[row sep=crcr]{%
1	1.26425117814664e-09\\
2	4.41694568937485e-09\\
3	1.11901624676687e-09\\
4	2.00779313640617e-09\\
5	2.92696712713986e-09\\
6	3.53780778066161e-09\\
7	1.77240078618648e-09\\
8	1.28417252499003e-09\\
9	5.88778158090964e-09\\
10	2.63315729446746e-09\\
11	2.30530005909832e-09\\
12	4.35581637783571e-08\\
13	1.30130520833394e-08\\
14	1.80494388648163e-09\\
15	2.48760627021141e-09\\
16	2.02269947210355e-09\\
17	9.70089435116095e-09\\
18	3.8962512187326e-08\\
19	1.54187771263982e-08\\
20	5.52408865668515e-09\\
21	5.08910103246949e-09\\
22	7.10493101920034e-09\\
23	2.45140846411628e-08\\
24	9.98584847150976e-09\\
25	1.04394265556703e-07\\
26	1.15316071329663e-09\\
27	1.30991966456835e-09\\
28	6.79151735870762e-09\\
29	3.45381233387086e-08\\
30	3.32512872560272e-09\\
31	4.54558170964736e-09\\
32	1.34895293356009e-09\\
33	5.85807091111922e-09\\
34	4.07223002105542e-09\\
35	1.07774175624214e-08\\
36	2.88526685584011e-08\\
37	1.18650273082861e-08\\
38	2.43613628589107e-09\\
39	4.35646038902812e-08\\
40	7.82085279544689e-09\\
41	2.50352753061352e-08\\
42	2.71788169979118e-09\\
43	8.85365184434468e-09\\
44	2.75738604272531e-08\\
45	3.52532837957068e-08\\
46	3.21548717679432e-08\\
47	1.46227430028402e-09\\
48	1.49810151979757e-08\\
49	2.42513412384201e-08\\
50	1.10154579190876e-07\\
51	1.71959433851873e-09\\
52	4.88404826467569e-09\\
53	2.78825064371469e-09\\
54	3.81164158052354e-09\\
55	1.08697370943425e-09\\
56	1.42896543372582e-09\\
57	9.11146248197058e-09\\
58	2.56414208240368e-08\\
59	1.42821308932445e-09\\
60	4.18337710640893e-08\\
61	2.51362597824991e-09\\
62	1.58726423861798e-09\\
63	2.54778245818193e-09\\
64	2.28442594933191e-09\\
65	1.76243109485581e-09\\
66	3.73254287207167e-08\\
67	5.24338963438041e-09\\
68	4.59286335039596e-09\\
69	1.15523415057906e-09\\
70	3.59821072025614e-09\\
71	2.52118257684296e-09\\
72	1.31255334646564e-07\\
73	1.05366431775145e-08\\
74	1.20918603022592e-08\\
75	7.5849905701352e-09\\
76	1.4940539581517e-08\\
77	3.28480802335184e-09\\
78	2.71188116992517e-09\\
79	1.78764686563986e-09\\
80	2.62466881469544e-08\\
81	1.21682813158397e-09\\
82	2.22673483219787e-09\\
83	3.84118179347359e-08\\
84	7.60745287502713e-09\\
85	7.10920469259536e-08\\
86	1.14018727735998e-07\\
87	2.53183535860425e-09\\
88	4.17301041039968e-09\\
89	2.28740106672019e-08\\
90	2.5466589112565e-09\\
91	4.11022580186901e-09\\
92	1.40713927519824e-09\\
93	9.50457646709851e-09\\
94	8.07630589162314e-09\\
95	2.8419140484531e-09\\
96	2.30941569707579e-09\\
97	2.69253871966659e-08\\
98	8.11458371493611e-09\\
99	2.92612675144054e-09\\
100	9.09059580879364e-09\\
};
\addlegendentry{Adaptive}
\addplot [color=mpigreen, dashdotted, very thick]
  table[row sep=crcr]{%
1	1e-05\\
2	1e-05\\
3	1e-05\\
4	1e-05\\
5	1e-05\\
6	1e-05\\
7	1e-05\\
8	1e-05\\
9	1e-05\\
10	1e-05\\
11	1e-05\\
12	1e-05\\
13	1e-05\\
14	1e-05\\
15	1e-05\\
16	1e-05\\
17	1e-05\\
18	1e-05\\
19	1e-05\\
20	1e-05\\
21	1e-05\\
22	1e-05\\
23	1e-05\\
24	1e-05\\
25	1e-05\\
26	1e-05\\
27	1e-05\\
28	1e-05\\
29	1e-05\\
30	1e-05\\
31	1e-05\\
32	1e-05\\
33	1e-05\\
34	1e-05\\
35	1e-05\\
36	1e-05\\
37	1e-05\\
38	1e-05\\
39	1e-05\\
40	1e-05\\
41	1e-05\\
42	1e-05\\
43	1e-05\\
44	1e-05\\
45	1e-05\\
46	1e-05\\
47	1e-05\\
48	1e-05\\
49	1e-05\\
50	1e-05\\
51	1e-05\\
52	1e-05\\
53	1e-05\\
54	1e-05\\
55	1e-05\\
56	1e-05\\
57	1e-05\\
58	1e-05\\
59	1e-05\\
60	1e-05\\
61	1e-05\\
62	1e-05\\
63	1e-05\\
64	1e-05\\
65	1e-05\\
66	1e-05\\
67	1e-05\\
68	1e-05\\
69	1e-05\\
70	1e-05\\
71	1e-05\\
72	1e-05\\
73	1e-05\\
74	1e-05\\
75	1e-05\\
76	1e-05\\
77	1e-05\\
78	1e-05\\
79	1e-05\\
80	1e-05\\
81	1e-05\\
82	1e-05\\
83	1e-05\\
84	1e-05\\
85	1e-05\\
86	1e-05\\
87	1e-05\\
88	1e-05\\
89	1e-05\\
90	1e-05\\
91	1e-05\\
92	1e-05\\
93	1e-05\\
94	1e-05\\
95	1e-05\\
96	1e-05\\
97	1e-05\\
98	1e-05\\
99	1e-05\\
100	1e-05\\
};
\addlegendentry{Tolerance}
\end{axis}
\end{tikzpicture}%
}
				\caption{Results for the Burgers' equation.}
				\label{fig:1}
			\end{figure}
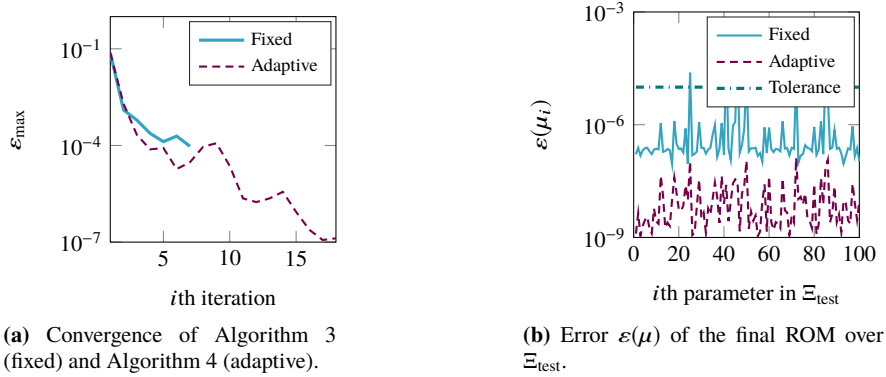	
			\setlength\fheight{1.7cm}
			\setlength\fwidth{6cm}			
				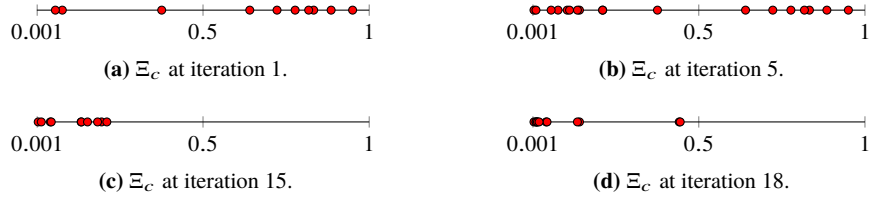
\begin{figure}[t]
					\centering
					\subfloat[\label{fig:2a}][$\Xi_c$ at iteration $1$.]{\begin{tikzpicture}
\pgfplotsset{ every non boxed x axis/.append style={x axis line style=-} }
\begin{axis}[
		width=\fwidth,
		height = \fheight,
		at={(0\fwidth,0\fheight)},
    xmin=0, xmax=1,
    axis x line=bottom,
		xtick={0.001,0.5, 1},
		xticklabel style={
						/pgf/number format/fixed,
						/pgf/number format/precision=4
		},
		scaled x ticks=false,
    hide y axis,    
    ymin=0,ymax=0.1,
    scatter/classes={%
        a={mark=o,draw=black}}
    ]

\addplot[scatter,only marks,
    mark size = 1.5pt,
    fill = red,
    scatter src=explicit symbolic]
table {
0.375682071744412 0
0.640664315302722 0
0.950066562825762 0
0.0766015274973595 0
0.777152729059402 0
0.832913016799297 0
0.0557509345421258 0
0.817859704793093 0
0.885466106614694 0
0.722624164296438 0
    };
\end{axis}

\end{tikzpicture}} \hfill
					\subfloat[\label{fig:2b}][$\Xi_c$ at iteration $5$.]{\begin{tikzpicture}
\pgfplotsset{ every non boxed x axis/.append style={x axis line style=-} }
\begin{axis}[
		width=\fwidth,
				height = \fheight,
		at={(0\fwidth,0\fheight)},
    xmin=0, xmax=1,
    axis x line=bottom,
		xtick={0.001,0.5, 1},
		xticklabel style={
						/pgf/number format/fixed,
						/pgf/number format/precision=4
		},
		scaled x ticks=false,		
    hide y axis,    
    ymin=0,ymax=0.1,
    scatter/classes={%
        a={mark=o,draw=black}}
    ]

\addplot[scatter,only marks,
    mark size = 1.5pt,
    fill = red,
    scatter src=explicit symbolic]
table {
0.375682071744412 0
0.640664315302722 0
0.950066562825762 0
0.0766015274973595 0 
0.777152729059402 0
0.832913016799297 0
0.0557509345421258 0
0.817859704793093 0
0.885466106614694 0
0.722624164296438 0
0.00413557918874456 0
0.0101207885434704 0
0.140961460432658 0
0.134884253020095 0
0.210890131530076 0
0.210579773735298 0
0.104129622118863 0
0.111132805121343 0
    };
\end{axis}

\end{tikzpicture}} \qquad
					\subfloat[\label{fig:2c}][$\Xi_c$ at iteration $15$.]{\begin{tikzpicture}
\pgfplotsset{ every non boxed x axis/.append style={x axis line style=-} }
\begin{axis}[
		width=\fwidth,
				height = \fheight,
		at={(0\fwidth,0\fheight)},
    xmin=0, xmax=1,
    axis x line=bottom,
		xtick={0.001,0.5, 1},
		xticklabel style={
						/pgf/number format/fixed,
						/pgf/number format/precision=4
		},
		scaled x ticks=false,		
    hide y axis,    
    ymin=0,ymax=0.1,
    scatter/classes={%
        a={mark=o,draw=black}}
    ]

\addplot[scatter,only marks,
    mark size = 1.5pt,
    fill = red,
    scatter src=explicit symbolic]
table {
0.00413557918874456 0
0.194757863995041 0
0.133208691040553 0
0.134071495392772 0
0.210579773735298 0
0.182714662890540 0
0.152480335207972 0
0.0133265555713069 0
0.0405482014859437 0
0.0433652936557696 0
    };
\end{axis}

\end{tikzpicture}}\hfill	
					\subfloat[\label{fig:2d}][$\Xi_c$ at iteration $18$.]{\begin{tikzpicture}
\pgfplotsset{ every non boxed x axis/.append style={x axis line style=-} }
\begin{axis}[
		width=\fwidth,
		height = \fheight,
		at={(0\fwidth,0\fheight)},
    xmin=0, xmax=1,
    axis x line=bottom,
		xtick={0.001,0.5, 1},
		xticklabel style={
						/pgf/number format/fixed,
						/pgf/number format/precision=4
		},
		scaled x ticks=false,		
    hide y axis,    
    ymin=0,ymax=0.1,
    scatter/classes={%
        a={mark=o,draw=black}}
    ]

\addplot[scatter,only marks,
    mark size = 1.5pt,
    fill = red,
    scatter src=explicit symbolic]
table {
0.00413557918874456 0
0.0133265555713069 0
0.0405482014859437 0
0.0433652936557696 0
0.0101207885434704 0
0.0108032502242695 0
0.0146348421601551 0
0.0194116538175655 0
0.140961460432658 0
0.134884253020095 0
0.441185870925416 0
0.444086729092199 0
    };
\end{axis}

\end{tikzpicture}}					 
					\caption{Burgers' equation: training set evolution.}
					\label{fig:2}
				\end{figure}
\subsection{Convection-diffusion equation}
	Next, we consider a $1$-D model of brain transport, originally discussed in \cite{BanK89} and also considered in \cite{morGre05, morZhaFLetal15}. The transport is modelled as a linear convection-diffusion PDE defined in the spatial domain $w \in \Omega \coloneqq [0\,,1]$ and for time $T \coloneqq [0\,, 1]$,
	\begin{align}
		\frac{dv}{dt} &= q_{1} \frac{\partial^{2} v}{\partial w^{2}} + q_{2} \frac{\partial v}{\partial w} - q_{2},
	\end{align}
where $v(w,t)$ is the state vector and the two parameters $(q_1, q_2) \in \mathscr{P} := [0.001\,,1] \times [0.5\,,5]$ are the diffusion and convection constants, respectively. The boundary conditions are given by,
	\begin{align*}
			v(w,0) &= \begin{cases}
							1, & w \leq 0.5\\
							0, & \text{otherwise}
					  \end{cases}, \qquad v(0,t) = v(1,t) = 0.
	\end{align*}
We discretize the equation using the FDM on a grid yielding $n = 800$. The output is calculated as the average value of the state in a small interval $\Omega_{o}$ centered around the midpoint of the domain at $w = 0.5$. 
	\begin{align*}
		y(t) \coloneqq \frac{1}{|\Omega_{o}|} \int_{\Omega_{o}} v(w,t) dw, \qquad \Omega_{o} \coloneqq [0.495\,,0.505].
	\end{align*}
				\begin{figure}[b!]
					\centering			
							\setlength\fheight{3cm}
							\setlength\fwidth{3cm}
							\subfloat[Convergence of~\cref{alg:2} (fixed) and~\cref{alg:3} (adaptive).]{\label{fig:3a}
%
%
\begin{tikzpicture}
\begin{axis}[%
width=\fwidth,
height=\fheight,
at={(0\fwidth,0\fheight)},
scale only axis,
xmin=0,
xmax=31, 
xlabel={$i$th iteration},
xlabel style={font=\color{white!15!black},font = \small},
ymode=log,
ymin=1e-07, 
ymax=100,
yminorticks=true,
ylabel={$\varepsilon_{\text{max}}$},
ylabel style={font=\color{white!15!black}, font = \small},
legend style={legend cell align=left, align=left, draw=white!15!black, font = \scriptsize},
axis background/.style={fill=white}
]
\addplot [color=mpiblue, thick]
  table[row sep=crcr]{%
1	2.18119137979877\\
2	0.608868261238195\\
3	0.157756626263737\\
4	0.250358318164483\\
5	0.0593845791638206\\
6	0.0404520881797726\\
7	0.0560816891937944\\
8	0.0150562853929223\\
9	0.0381230451755694\\
10	0.0117987380612721\\
11	0.0046517024279732\\
12	0.00273125153811437\\
13	0.00516806390803082\\
14	0.0086998888050391\\
15	0.00678518943841794\\
16	0.00678695838041653\\
17	0.00665736156356646\\
};
\addlegendentry{Fixed}

\addplot [color=mpired, thick,densely dashed]
  table[row sep=crcr]{%
1	2.18119137979877\\
2	0.608868261238195\\
3	0.418382688477269\\
4	0.173720330218498\\
5	0.0828898009761649\\
6	0.053751824648757\\
7	0.00925630295294672\\
8	0.00152726011976441\\
9	0.00114256420724411\\
10	0.00107652852076716\\
11	0.000862474203860055\\
12	0.00178815919496108\\
13	0.00157243223011358\\
14	0.00102497177482498\\
15	0.00093933112682543\\
16	0.0012852859618206\\
17	0.000157519132988445\\
18	0.000150788375443509\\
19	0.000102415769541912\\
20	5.53382214611148e-05\\
21	2.74349616110172e-05\\
22	2.80952493550083e-05\\
23	2.73056280922798e-05\\
24	2.30063714754412e-05\\
25	1.57207765281642e-05\\
26	6.05160246424767e-06\\
27	1.17396792156077e-05\\
28	1.10198059620313e-05\\
29	9.94698328098464e-06\\
30	2.20544922980534e-06\\
31	1.52235203369395e-06\\
};
\addlegendentry{Adaptive}
\end{axis}
\end{tikzpicture}%
}	\hfill
							\subfloat[Error $\varepsilon(\mu)$ of the final ROM over $\Xi_{\text{test}}$.]{\label{fig:3b}\input{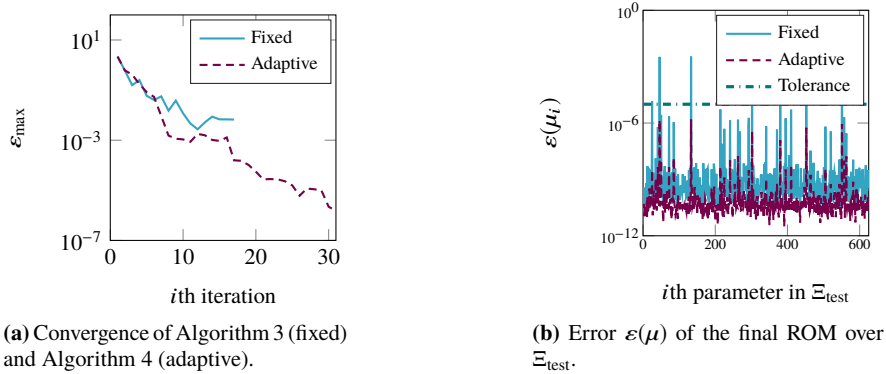}}			
							\caption{Results for the Convection-diffusion equation.}
							\label{fig:3} 
				\end{figure}
				\setlength\fheight{2.5cm}
				\setlength\fwidth{2.5cm}
				\begin{figure}[b!]
				\subfloat[\label{fig:cd-a}][$\Xi_c$ at iteration $1$.]{
%
%
\begin{tikzpicture}
\pgfplotsset{every x tick label/.append style={font=\small, yshift=-0.5ex}}
\begin{axis}[%
width=\fwidth,
height=\fheight,
at={(0\fwidth,0\fheight)},
scale only axis,
xmin=0,
xmax=1,
xlabel style={font=\color{white!15!black}, font = \small},
xlabel={$q_{1}$},
xtick={0.001, 0.5, 1},
xlabel near ticks,
		xticklabel style={
						/pgf/number format/fixed,
						/pgf/number format/precision=4
		},
		scaled x ticks=false,
ymin=0.5,
ymax=5,
ylabel style={font=\color{white!15!black}, font = \small},
ylabel={$q_{2}$},
ylabel near ticks,
ytick={0.5, 3, 5},
axis background/.style={fill=white},
]
\addplot[only marks, mark=asterisk, mark options={},thick, mark size=2pt, draw=mpiblue] table[row sep=crcr]{%
x	y\\
0.375682071744412	1.80531376561465\\
0.640664315302722	0.702536775475013\\
0.950066562825762	1.99484490798473\\
0.0766015274973595	3.81439947308247\\
0.777152729059402	4.27645168344848\\
0.832913016799297	3.69212726575844\\
0.0557509345421258	1.34780502474865\\
0.817859704793093	4.93550851960727\\
0.885466106614694	2.2460148862045\\
0.722624164296438	2.76419516865737\\
0.00355347569162105	3.9358864638591\\
0.981218078443031	3.57679602275169\\
0.344076430916395	3.74880821631083\\
0.0956651282880281	1.4636844084375\\
0.39524794610614	3.93391660538238\\
0.00593997988690899	0.838472699216402\\
0.736970123437002	3.24900815767537\\
0.955851062616075	2.41533691095288\\
0.820788761943263	1.55457989151737\\
0.345638013090133	2.97062922911623\\
0.379165913347508	3.09604888616082\\
0.784978851354912	2.96607316268804\\
0.0871452815243885	1.19710174808105\\
0.546527863863037	3.34652666474077\\
0.163063794268542	0.99114537844222\\
};
\end{axis}
\end{tikzpicture}%
}
					\hfill
					\subfloat[\label{fig:cd-c}][$\Xi_c$ at iteration $20$.]{
%
%
\begin{tikzpicture}
\pgfplotsset{every x tick label/.append style={font=\small, yshift=-0.5ex}}
\begin{axis}[%
width=\fwidth,
height=\fheight,
at={(0\fwidth,0\fheight)},
scale only axis,
xmin=0,
xmax=1,
xlabel style={font=\color{white!15!black}, font = \small},
xlabel={$q_{1}$},
xlabel near ticks,
		xticklabel style={
						/pgf/number format/fixed,
						/pgf/number format/precision=4
		},
		scaled x ticks=false,
xtick={0.001, 0.5, 1},
ymin=0.5,
ymax=5,
ylabel style={font=\color{white!15!black}, font = \small},
ylabel={$q_{2}$},
ylabel near ticks,
ytick={0.5, 3, 5},
axis background/.style={fill=white},
]
\addplot[only marks, mark=asterisk, mark options={}, mark size=2pt, draw=brown] table[row sep=crcr]{%
x	y\\
0.640664315302722	0.702536775475013\\
0.950066562825762	1.99484490798473\\
0.0557509345421258	1.34780502474865\\
0.885466106614694	2.2460148862045\\
0.722624164296438	2.76419516865737\\
0.981218078443031	3.57679602275169\\
0.344076430916395	3.74880821631083\\
0.0956651282880281	1.4636844084375\\
0.39524794610614	3.93391660538238\\
0.00593997988690899	0.838472699216402\\
0.955851062616075	2.41533691095288\\
0.820788761943263	1.55457989151737\\
0.784978851354912	2.96607316268804\\
0.0871452815243885	1.19710174808105\\
0.163063794268542	0.99114537844222\\
0.103461538461538	3.73076923076923\\
0.103461538461538	3.84615384615385\\
0.0266153846153846	2.92307692307692\\
0.0778461538461538	2.92307692307692\\
0.0778461538461538	3.03846153846154\\
0.0522307692307692	2.92307692307692\\
0.0266153846153846	3.03846153846154\\
0.487692307692308	4.42307692307692\\
0.487692307692308	4.30769230769231\\
0.103461538461538	4.19230769230769\\
0.0778461538461538	4.30769230769231\\
0.129076923076923	5\\
0.154692307692308	5\\
0.001	0.615384615384615\\
0.001	0.5\\
0.0266153846153846	0.5\\
0.0266153846153846	0.615384615384615\\
0.001	0.730769230769231\\
0.001	0.846153846153846\\
0.001	1.42307692307692\\
0.001	1.30769230769231\\
0.103461538461538	0.615384615384615\\
0.129076923076923	0.615384615384615\\
};
\end{axis}
\end{tikzpicture}%
}		
					\hfill
					\subfloat[\label{fig:cd-d}][$\Xi_c$ at iteration $31$.]{
%
%
\definecolor{mycolor1}{rgb}{0.00000,1.00000,1.00000}%
\begin{tikzpicture}
\pgfplotsset{every x tick label/.append style={font=\small, yshift=-0.5ex}}
\begin{axis}[%
width=\fwidth,
height=\fheight,
at={(0\fwidth,0\fheight)},
scale only axis,
xmin=0,
xmax=1,
xlabel style={font=\color{white!15!black}, font = \small},
xlabel={$q_{1}$},
xtick={0.001, 0.5, 1},
xlabel near ticks,
		xticklabel style={
						/pgf/number format/fixed,
						/pgf/number format/precision=4
		},
		scaled x ticks=false,
ymin=0.5,
ymax=5,
ylabel style={font=\color{white!15!black}, font = \small},
ylabel={$q_{2}$},
ylabel near ticks,
ytick={0.5, 3, 5},
axis background/.style={fill=white},
]
\addplot[only marks, mark= asterisk, mark options={}, mark size=2pt, draw= mpired] table[row sep=crcr]{%
x	y\\
0.001	0.615384615384615\\
0.001	0.5\\
0.0266153846153846	0.5\\
0.001	0.730769230769231\\
0.001	0.846153846153846\\
0.001	1.42307692307692\\
0.001	1.30769230769231\\
0.103461538461538	0.615384615384615\\
0.001	5\\
0.001	2\\
0.0266153846153846	2\\
0.001	1.65384615384615\\
0.0266153846153846	2.23076923076923\\
0.001	1.07692307692308\\
0.001	0.961538461538462\\
0.001	3.38461538461538\\
0.0522307692307692	2.11538461538462\\
0.001	2.23076923076923\\
0.0522307692307692	0.5\\
0.257153846153846	0.5\\
0.282769230769231	0.5\\
1	5\\
1	4.88461538461539\\
0.0778461538461538	0.5\\
0.0266153846153846	5\\
};
\end{axis}
\end{tikzpicture}%
}					 
					\caption{Convection-diffusion equation: training set evolution.}
					\label{fig:4}
				\end{figure}
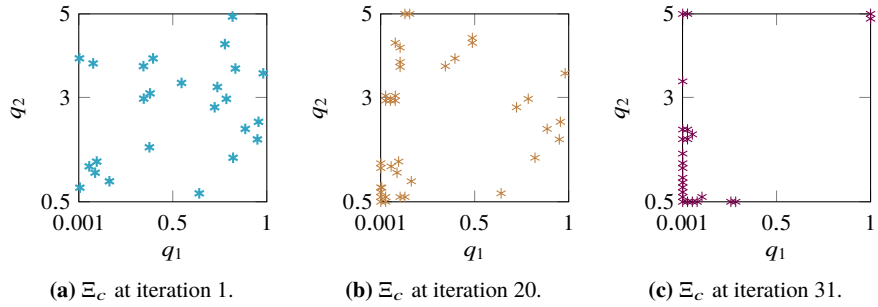
The tolerance is set as $\texttt{tol} = 10^{-5}$. For~\cref{alg:2}, we set $\Xi$ as a set including $25$ random samples from the parameter domain $\mathscr{P}$, picked using the \texttt{rng} command and making use of \texttt{twister}, with a seed of $112$. The fine training set $\Xi_f$ for~\cref{alg:3} includes $1600$ equidistant samples in $\mathscr{P}$. We choose a test set $\Xi_{\text{test}}$ consisting of $625$ random samples. The test parameters are generated using the same random number  generator as for the Burgers' equation example, viz., \texttt{simdTwister} with a seed of $200$.
The error surrogate model $s(\mu)$ in~\cref{alg:3} is constructed using IMQ kernel function, where cross validation LOOCV has been applied to specify the shape parameter. 

\cref{fig:3a} shows the maximal error decay ($\varepsilon_{\text{max}}$) over $\Xi_{\text{test}}$, computed at each iteration of~\cref{alg:2} and~\cref{alg:3}, respectively. 
It is clear that adaptively enriching the training set leads to orders of magnitude faster convergence to the required accuracy.
For \cref{alg:2}, the ROM size is 38 and for \cref{alg:3}, it is 63.
\cref{fig:3b} plots the error $\varepsilon(\mu)$ (\cref{eq:testerrorplotfinal}) of the final ROM at every parameter in the test set $\Xi_{\text{test}}$, for the case of an adaptive training set and a fixed training set. We see that for some test parameters, the required tolerance is not met by the ROM computed using the fixed training set. In \cref{fig:4}, the evolution of the training set is shown at different stages of~\cref{alg:3}. It can be seen that samples from the left boundary of the parameter domain are added or kept. This has physical sense as this corresponds to the lower viscosity regions where the convective part of the solution dominates and the solution is `less smooth'. Thus, more basis functions are required for a good approximation.

In order to achieve a better ROM by using the fixed training set, we tried a $\Xi$ with 225 random samples for~\cref{alg:2}, leading to a ROM with error below the tolerance. In~\cref{tab:convdiff_runtime}, we provide the runtime results of both algorithms, where~\cref{alg:2} uses the refined training set with 225 samples. On the other hand, the initial coarse training set and the fine training set for~\cref{alg:3} remain the same. ~\cref{alg:3} with adaptive sampling outperforms~\cref{alg:2}, though both approaches produce accurate ROMs.  
\subsection{A thermal model}
The final example is the Thermal model taken from the $\texttt{MORWiki}$ benchmark collection\footnote{\texttt{https://morwiki.mpi-magdeburg.mpg.de/morwiki/index.php/Thermal\_Model}}. It models the heat transfer in a microthruster unit. The system is governed by the following semi-discretized ODE,
		\begin{align*}
			\mathbf{E} \dot{x} = (\mathbf{A} - \sum_{i=1}^{3} h_{i} \mathbf{A}_{i})x + \mathbf{B} u.
		\end{align*}
The model has been discretized in space using the FEM and the matrices are available from the Oberwolfach Benchmark Collection hosted at \texttt{MORWiki}. The dimension of the system is $n = 4257$. The parameters $\{ h_{i} \}_{i=1}^{3}$	are the film coefficients at the top, bottom and side of the microthruster unit, respectively.
They have a range between $[1\,,10^{8}]$ in our simulations. The output is the temperature at the center of the polysilicon heater in the unit, which corresponds to the first row of the output matrix $C$.

\begin{figure}[t]
	\centering		
			\setlength\fheight{3cm}
			\setlength\fwidth{3cm}
			\subfloat[Convergence of~\cref{alg:2} (fixed) and~\cref{alg:3} (adaptive).]{\label{fig:5a}
%
%
\definecolor{mycolor1}{rgb}{0.00000,0.44700,0.74100}%
\begin{tikzpicture}

\begin{axis}[%
width=\fwidth,
height=\fheight,
at={(0\fwidth,0\fheight)},
scale only axis,
xmin=0,
xmax=24,
xlabel style={font=\color{white!15!black},font = \small},
xlabel={$i$th iteration},
ymode=log,
ymin=0.0001,
ymax=100,
yminorticks=true,
ylabel style={font=\color{white!15!black}, font = \small},
ylabel={$\varepsilon_{\text{max}}$},
axis background/.style={fill=white},
legend style={legend cell align=left, align=left, draw=white!15!black, font = \scriptsize}
]
\addplot [color=mpiblue, thick]
  table[row sep=crcr]{%
1	37.8073991211035\\
2	4.5927752410146\\
3	1.08194748980317\\
4	1.08215685131126\\
5	0.349569896188333\\
6	0.0445088624887898\\
7	0.0443549638064891\\
8	0.0443196410233602\\
9	0.0442674910077981\\
10	0.0440904339590043\\
11	0.0440378250205333\\
12	0.0440374822252663\\
13	0.0440369922090192\\
14	0.0440210537481677\\
15	0.0439338003728768\\
16	0.04384726408663\\
17	0.0438284040668422\\
18	0.0436106685578472\\
19	0.012014794345297\\
20	0.00156747663604592\\
21	0.00156776527168344\\
22	0.00156763442656467\\
23	0.00156763271823205\\
};
\addlegendentry{Fixed}
\addplot [color=mpired, thick,densely dashed]
table[row sep=crcr]{%
	1	60.3640336123813\\
	2	47.8347373270495\\
	3	47.738100254001\\
	4	47.6728070714911\\
	5	47.262006706903\\
	6	46.6209979190177\\
	7	41.4345236481562\\
	8	0.114736610425806\\
	9	0.0391838990838522\\
	10	0.039060105203028\\
	11	0.00847982958750539\\
	12	0.00450958804384907\\
	13	0.0043932524503017\\
	14	0.00230879300763227\\
	15	0.00192494658462178\\
	16	0.00192151770679803\\
	17	0.00181148197959449\\
	18	0.000392455253364701\\
	19	0.000384702967670239\\
	20	0.0011103168285809\\
	21	0.00100243133604299\\
	22	0.000583235170340634\\
	23	0.000599071909504341\\
	24	0.00057748765736975\\
};
\addlegendentry{Adaptive}
\end{axis}
\end{tikzpicture}%
} \hfill
			\subfloat[Error $\varepsilon(\mu_i)$ of the final ROM over $\Xi_{\text{test}}$.]{\label{fig:5b}		\input{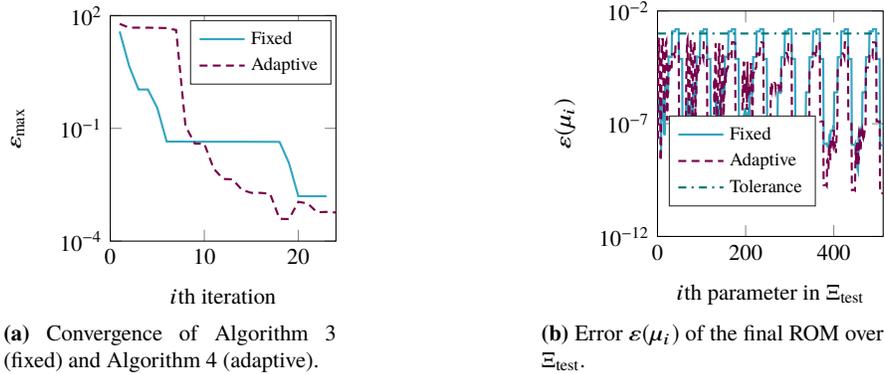}} 			
\caption{Results for the thermal model.}
\label{fig:5}
\end{figure}

Since this is a three-parameter problem with very large parameter domain, we use a large fixed training set ($\Xi$) with $6^3$ logarithmically equidistant parameter samples for~\cref{alg:2} for a fair comparison. For \cref{alg:3}, the initial coarse training set includes only $N_c = 10$ randomly chosen samples. To generate these random samples, we first generate a uniform sampling of each parameter $h_j, j=1,2,3$, given as, $h_{ij} = 10^{\frac{i}{N_{c}/8}},\, i = 1, 2, \ldots, N_{c}$, which leads to a matrix $(h_{ij})\in \mathbb R^{N_c\times 3}$.
Then, the rows of the matrix are permuted using the \texttt{rng} and \texttt{randperm} commands. The seeds and random number generators used are $100\,  (\texttt{twister}), 120\, (\texttt{combRecursive})\,\, \text{and}\,\,600\, (\texttt{combRecursive})$, respectively. In this way, we get random samples which spread through the 3-D parameter domain. The fine training set $\Xi_f$ for \cref{alg:3} is made up of $16^3$ equidistant samples, while the test set $\Xi_{\text{test}}$ consists of  $8^3$ logarithmically equidistant samples. For the RBF interpolation, Thin-plate spline kernel is used. No shape parameter needs to be determined. The tolerance for this model is set as $\texttt{tol} = 10^{-3}$.
In \cref{fig:5a}, we show error decay $\varepsilon_{\text{max}}$ over $\Xi_{\text{test}}$ at each iteration of~\cref{alg:2} and~\cref{alg:3}, respectively. \cref{alg:2} takes 23 iterations to converge with the fixed training set $\Xi$. The resulting order of the ROM is $44$. However, the maximum error $\varepsilon_{\text{max}}$ (\cref{eq:errorplot}) over the test set $\Xi_{\text{test}}$ is still above the tolerance. \cref{alg:3} converges in 24 iterations, to a ROM of order $74$. The maximum error $\varepsilon_{\text{max}}$ over $\Xi_{\text{test}}$ is below the tolerance upon convergence. We note that, in the first few iterations, \cref{alg:2} has a faster convergence in comparison to \cref{alg:3}. However, since the training set is fixed, the convergence eventually saturates. In \cref{fig:5b}, we plot the error $\varepsilon(\mu)$ of the final ROM over 
$\Xi_{\text{test}}$, computed by the two algorithms respectively. \cref{alg:3} once again outperforms \cref{alg:2}.

In \cref{tab:thermal_runtime}, we provide the runtime comparison between \cref{alg:2} and \cref{alg:3} for the thermal model, where an obvious speed-up by~\cref{alg:3} is observed. More importantly, with the reduced runtime,~\cref{alg:3} produces a ROM with sufficient accuracy, whereas, the ROM computed by~\cref{alg:2} still does not meet the accuracy requirement. This further justifies the motivation of using adaptive sampling for models with 
two or more parameters, especially when the parameter domain is very large. 

\begin{table}[t]
	\centering		
	\subfloat[Convection-diffusion example.]{\label{tab:convdiff_runtime} \scriptsize
	\begin{tabular}{|c|c|}
		\hline
		Algorithm\Tstrut & Runtime (seconds)\Tstrut \\[0.5ex]
		\hline
		\hline
		Fixed training set (\cref{alg:2})\Tstrut & $74.09$\Tstrut \\[0.5ex] 
		Adaptive training set (\cref{alg:3}) & $41.52$\\[0.5ex]	 
		\hline
	\end{tabular}}	\hfill
	\subfloat[Thermal example.]{\label{tab:thermal_runtime} \scriptsize
	\begin{tabular}{|c|c|}
		\hline
		Algorithm\Tstrut & Runtime (seconds)\Tstrut \\[0.5ex]
		\hline
		\hline
		Fixed training set (\cref{alg:2})\Tstrut & $151.16$\Tstrut \\[0.5ex] 
		Adaptive training set (\cref{alg:3}) & $38.76$\\[0.5ex]	 
		\hline
	\end{tabular}}
	\caption{Runtime comparison between \cref{alg:2} and \cref{alg:3}.}
	\label{tab:runtime}       
\end{table}

\section{Conclusion}
\label{sec:5}
Using a fixed training set in the POD-Greedy algorithm may either entail high computational costs or lead to large errors if the parameter domain is not properly sampled. In this work, we introduce a method to construct the training set for the RB method in an adaptive manner. To achieve this, we make use of an error surrogate model based on RBF interpolation. 
	
The use of an error surrogate model enables sufficient exploration of the parameter space, since it avoids the need to solve the ROM at every parameter in the training set, as required by the error estimator. Furthermore, we implement a cross validation strategy in order to choose good shape parameters for the RBF kernels.
The proposed algorithm is tested on several examples and it is shown to be effective in constructing a ROM with required accuracy. Furthermore, the adaptive parameter sampling method is integrated with the adaptive POD-Greedy-(D)EIM algorithm in~\cite{morCheFB19a} to achieve a fully adaptive scheme for the RB method.
	
	As future work, we propose to extend the adaptive sampling approach to the frequency domain model reduction methods, such as multi-moment matching \cite{morFenB14}, in order to adaptively sample the interpolation points.
	
 \bibliographystyle{spmpsci}
 \bibliography{myref}
\end{document}